%% file: Fast_Sinkhorn_II.tex
\documentclass[review,hidelinks,onefignum,onetabnum]{siamart220329}

\usepackage{setspace}
\input{ex_shared}

\newcommand{\brac}[1]{\left(#1\right)}
\newcommand{\norm}[1]{\left\Vert#1\right\Vert}

\newcommand{\bx}{\boldsymbol{x}}
\newcommand{\by}{\boldsymbol{y}}

\newcommand{\bp}{\boldsymbol{p}}
\newcommand{\bq}{\boldsymbol{q}}
\newcommand{\br}{\boldsymbol{r}}
\newcommand{\bs}{\boldsymbol{s}}
\newcommand{\mbr}{\mathbb{R}}
\newcommand{\mbn}{\mathbb{N}}
\newcommand{\bu}{\boldsymbol{u}}
\newcommand{\bv}{\boldsymbol{v}}
\newcommand{\bgamma}{\boldsymbol{\gamma}}
\newcommand{\balpha}{\boldsymbol{\alpha}}
\newcommand{\bbeta}{\boldsymbol{\beta}}
\newcommand{\bphi}{\boldsymbol{\phi}}
\newcommand{\bpsi}{\boldsymbol{\psi}}
\newcommand{\veps}{\varepsilon}
\ifpdf
\hypersetup{
  pdftitle={An Example Article},
  pdfauthor={D. Doe, P. T. Frank, and J. E. Smith}
}
\fi

\begin{document}
\begin{sloppypar}
\maketitle

\begin{abstract}
    In our previous work [arXiv:2202.10042], the complexity of Sinkhorn iteration is reduced from $O(N^2)$ to the optimal $O(N)$ by leveraging the special structure of the kernel matrix. In this paper, we explore  the special structure of kernel matrices by defining and utilizing the properties of the \textbf{L}ower-\textbf{Co}l\textbf{L}inear \textbf{T}riangular Matrix (L-CoLT matrix) and \textbf{U}pper-\textbf{Co}l\textbf{L}inear \textbf{T}riangular Matrix (U-CoLT matrix). We prove that (1) L/U-CoLT matrix-vector multiplications can be carried out in $O(N)$ operations; (2) both families of matrices are closed under the Hadamard product and matrix scaling. These properties help to alleviate  two key difficulties for reducing the complexity of the Inexact Proximal point method (IPOT), and allow us to  significantly reduce the number of iterations  to $O(N)$. This yields the Fast Sinkhorn II (FS-2) algorithm for accurate computation of optimal transport with low algorithm complexity and fast convergence. Numerical experiments are presented to show the effectiveness and efficiency of our approach.
\end{abstract}

\begin{keywords}
	Optimal Transport, Wasserstein-1 metric, Sinkhorn algorithm, IPOT method, FS-2 algorithm
\end{keywords}

\begin{MSCcodes}
    49M25; 65K10
\end{MSCcodes}

\section{Introduction}
\label{sec:1}
The Wasserstein metric, broadly used in  optimal transport theory with applications in many fields including machine learning, quantifies the dissimilarity between two probabilistic distributions. Many methods have been proposed to compute the Wasserstein metrics directly, such as the linear programming methods~\cite{pele2009fast,li2020asymptotically,yang2021fast}, combinatorial methods~\cite{santambrogio2015optimal}, solving the Monge-Amp\`ere equations~\cite{froese2011convergent,froese2012numerical,benamou2014numerical}, via Benamou-Brenier formulation \cite{benamou2000computational, li2018parallel} and the proximal splitting methods~\cite{combettes2011proximal,metivier2016optimal}. In recent years, several approximation techniques in optimal transport for high-dimensional distributions have also been proposed \cite{meng2019large,meng2020sufficient}.

The Sinkhorn algorithm~\cite{cuturi2013sinkhorn,sinkhorn1967diagonal} is a popular $O(N^2)$ algorithm to approximate the Wasserstein metric~\cite{peyre2019computational} by minimizing the entropy regularized optimal transport (OT) problem. In~\cite{liao2022fast}, the FS-1 algorithm is proposed to solve entropy regularized OT in $O(N)$ time by leveraging the special structure of the Sinkhorn kernel matrix of the Wasserstein-1 metric. The solution of entropy regularized OT approximates the accurate OT solution only if the regularization parameter is sufficiently small. However, small regularization parameters lead to numerical instability and excessive iterations~\cite{franklin1989scaling}. This causes the slow convergence of the Sinkhorn algorithm.

The Inexact Proximal point method~\cite{xie2020fast} for the Optimal Transport problem (IPOT) has been proposed to address this challenge. It regularizes the original OT by introducing the proximal point term and solves a series of successive subproblems. Only fairly mild regularization parameters are required to ensure the method's fast convergence to the accurate OT solution in an $O(N^2)$ algorithm. The goal of this paper is to construct a new method to accurately compute OT solutions with good convergence behavior and $O(N)$ algorithm complexity by combing the IPOT method and the FS-1 algorithm. Note the two key steps in the IPOT method make it hard to reduce the complexity to $O(N)$: the matrix Hadamard product (\cref{alg:IPOT}, line 4) and the matrix scaling (\cref{alg:IPOT}, line 8). For general matrices, the complexity of the above operations is both $O(N^2)$. Moreover, these operations may destroy the special structure of the kernel matrix \cite{liao2022fast}, making it impossible for us to implement matrix-vector multiplication with $O(N)$ cost.

 We will explore the special structure of kernel matrices by defining and exploiting the properties of the \textbf{L}ower-\textbf{Co}l\textbf{L}inear \textbf{T}riangular Matrix (L-CoLT matrix) and \textbf{U}pper-\textbf{Co}l\textbf{L}inear \textbf{T}riangular Matrix (U-CoLT matrix). For these matrices, we can realize the matrix-vector multiplication with $O(N)$ cost by using the idea of dynamic programming similar to \cite{liao2022fast}. Next, we show that each L/U-CoLT matrix can be represented by two vectors of dimension $N$. Furthermore, we prove the closure of families of L/U-CoLT matrices to matrix Hadamard product and matrix scaling. This means that the special structure of the kernel matrix is preserved by matrix Hadamard product and matrix scaling, so we can still implement matrix-vector multiplication (\cref{alg:IPOT}, lines 6-7) with $O(N)$ cost. On the other hand, by updating two representation vectors of the L/U-CoLT matrix, we can also implement matrix Hadamard product (\cref{alg:IPOT}, line 4) and matrix scaling (\cref{alg:IPOT}, line 8) with $O(N)$ cost. Consequently, the Fast Sinkhorn II (FS-2) algorithm is developed, which integrates the advantages of both IPOT and FS-1. Moreover, we also find that the FS-2 algorithm has the advantage in reducing the space complexity since all the matrices can be represented by vectors. Due to these benefits, one can expect that our FS-2 could be applied in various fields, e.g., machine learning~\cite{meng2019large,meng2020sufficient,goodfellow2014generative,lin2021wasserstein}, image processing~\cite{rubner2000earth,museyko2009application}, inverse problems~\cite{chen2018quadratic,engquist2020quadratic,yang2018application,heaton2020wasserstein}, density function theory\cite{hu2021global,buttazzo2012optimal,cotar2013density}.

The rest of the paper is organized as follows. In \cref{sec:2} , the basics of the Wasserstein-1 metric and the IPOT method are briefly reviewed. After presenting the definition, properties, and fast matrix-vector multiplications of the L/U-CoLT matrix in \cref{sec:ctm}, we apply them to accelerate the IPOT method, thus developing the FS-2 algorithm in \cref{sec:fs2}. In \cref{sec:fs2hd}, the FS-2 algorithm is extended to high dimensions. The numerical experiments are performed to verify our conclusions in \cref{sec:5}. We conclude the paper in \cref{sec:6}.

\section{The Wasserstein-1 metric and the IPOT method}
\label{sec:2}
Given two unit discrete distributions $\bu$ and $\bv$,
\begin{equation*}
	\bu = (u_1,u_2,\cdots u_N)^\top\in\mbr^{N}, \quad  \bv = (v_1,v_2,\cdots, v_N)^\top \in\mbr^{N},
\end{equation*}
where $u_i \geq 0$, $v_j \geq 0$, and $\sum_i u_i = \sum_j v_j = 1$.
The Wasserstein-1 distance between them is defined as~\cite{peyre2019computational} 
\begin{equation}
    \label{eqn:ot}
    W_1\brac{\bu,\bv}=\min_{\Gamma\mathbf{1}=\bu,\Gamma^T\mathbf{1}=\bv}\left\langle C,\;\Gamma\right\rangle,
\end{equation}
where $C=[c_{ij}]\in\mbr^{N\times N}$ is the cost matrix. The element $c_{ij}=\Vert \bx_i-\by_j\Vert_1$ represents the cost of transporting the unit mass from position $\bx_i$ to position $\by_j$ and the variable $\Gamma=[\gamma_{ij}]\in\mbr^{N\times N}$ to be optimized is the transport plan. Here, the Frobenius inner product $\left\langle A,B\right\rangle=\sum_{i,j}a_{ij}b_{ij}$, where $A=[a_{ij}],B=[b_{ij}]$ are real-valued matrices.

The Sinkhorn algorithm~\cite{cuturi2013sinkhorn,sinkhorn1967diagonal} solves an entropy regularized OT problem to obtain an approximate result of \eqref{eqn:ot}. However, the small regular parameter required by good approximation leads to a slow convergence rate and numerical instability. To avoid this problem, the proximal point iteration \eqref{eqn:ipotproblem} is developed to solve \eqref{eqn:ot} accurately ~\cite{xie2020fast}. It begins with a transport map $\Gamma^{(0)}$ and iteratively solves the following minimization problem
\begin{equation}
    \label{eqn:ipotproblem}
    \Gamma^{\brac{t+1}}=\mathop{\arg\min}\limits_{\Gamma\mathbf{1}=\bu,\Gamma^T\mathbf{1}=\bv}\left\langle C,\;\Gamma\right\rangle+\delta^{\brac{t}}D_h\brac{\Gamma,\;\Gamma^{\brac{t}}},
\end{equation}
where $D_h$ is Bregman divergence,  taken in the form of the KL divergence in \cite{xie2020fast}, 
\[
D_h\brac{A,\;B}=\sum_{i,j}\brac{a_{ij}\ln\frac{a_{ij}}{b_{ij}}-a_{ij}+b_{ij}}
\]
and $\delta^{\brac{t}}$ is the regular parameter. The Lagrangian of the above equation writes
    \begin{equation*}\label{Laplace}
        L(\Gamma,\balpha,\bbeta) = \left\langle C,\;\Gamma\right\rangle+\delta^{\brac{t}}D_h\brac{\Gamma,\;\Gamma^{\brac{t}}}+\balpha^T\brac{\Gamma\mathbf{1}-\bu}+\bbeta^T\brac{\Gamma^T\mathbf{1}-\bv}.
    \end{equation*} 
    Taking derivative of the Lagrangian with respect of $\gamma_{ij}$ directly leads to
    \begin{equation*}
        \gamma_{ij} = 
        e^{-\alpha_i/\delta^{\brac{t}}} 
        Q_{ij}^{(t)}
        e^{-\beta_j/\delta^{\brac{t}}}, \;\;
        \textrm{where} \;\; Q_{ij}^{(t)} =\gamma^{\brac{t}}_{ij}e^{-c_{ij}/\delta^{\brac{t}}} > 0.
    \end{equation*}
    Denoting $\odot$ as the Hadamard product, $Q^{(t)}=K\odot\Gamma^{(t)}$ and $K=[e^{-c_{ij}/\delta^{(t)}}]\in\mbr^{N\times N}$ is the kernel matrix.
    Letting  $\phi_i=e^{-\alpha_i/\delta^{\brac{t}}}$,
    $\psi_j=e^{-\beta_j/\delta^{(t)}}$, and vectors $\bphi=(\phi_i)$ and $\bpsi=(\psi_j)$, one obtains
    \begin{equation}\label{solution}
        {\rm diag}(\bphi) Q^{(t)} {\rm diag}(\bpsi) \mathbf{1}  = \bu ,\quad {\rm diag}(\bpsi) Q^{(t)\top} {\rm diag}(\bphi) \mathbf{1} = \bv.
    \end{equation}
    By iteratively updating vectors $\bphi$ and $\bpsi$
    \begin{equation}\label{update_sinkhorn}
        \bpsi^{(t,\ell+1)} = \bv\oslash (Q^{(t)\top}\bphi^{(t,\ell)}),\quad
        \bphi^{(t,\ell+1)} = \bu\oslash (Q^{(t)}\bpsi^{(t,\ell+1)}),
    \end{equation}
    one can obtain an accurate solution for the original OT problem \eqref{eqn:ot}. Here $\oslash$ represents pointwise division, $t$ is the proximal iteration step (outer iteration) and $\ell$ it the Sinkhorn-type iteration step (inner iteration). The pseudo-code of IPOT is shown in \cref{alg:IPOT}.

\begin{algorithm} 
\setstretch{1.2} 
    \caption{IPOT}
    \label{alg:IPOT}
    \hspace*{0.02in} {\bf Input:} $\bu,\;\bv\in\mbr^{N}$;\ $K = e^{-C/\delta}\in\mbr^{N\times N}$;\ $L,\mathrm{itr\_max}\in\mbn^+$ \\
    \hspace*{0.02in} {\bf Output:} $W_1(\bu,\bv)$
    \begin{algorithmic}[1]
			\State $\bphi,\;\bpsi\gets \frac{1}{N}\mathbf{1}_N$
			\State $\Gamma=\mathbf{1}_N\mathbf{1}^T_N$
			\For{$t=1\;:\;\mathrm{itr\_max}$}
			\State $Q\gets K\odot\Gamma$
			\For{$\ell=1\;:\;L$}
			\State $\bpsi\gets\bv\oslash\brac{Q^T\bphi}$
			\State $\bphi\gets\bu\oslash\brac{Q\bpsi}$
			\EndFor
			\State $\Gamma\gets\text{diag}\brac{\bphi}Q\text{diag}\brac{\bpsi}$
			\EndFor
			\noindent\Return $W_1(\bu,\bv)$
		\end{algorithmic}
\end{algorithm}

\section{The Collinear Triangular Matrix}
\label{sec:ctm}

\subsection{Definition and Fast Matrix-Vector Multiplication}

\begin{definition}[Lower/Upper-Collinear Triangular Matrix]
\label{def:colt}
A lower triangular matrix is called a {\bf Lower-Collinear Triangular Matrix}(L-CoLT matrix) if its corresponding entries on any two rows (columns) have  the same--column (row) independent-- ratio except those dividing by $0$. Specifically, the N-dimensional L-CoLT matrix set is defined as follows: 
\begin{equation}\label{eqn:coltdef}
    \mathcal{C}_L^N\!=\!\left\{M\in\mbr^{N\times N}\;\vert\; m_{\!i+1,j\!}/m_{\!i,j\!}=r_i,\; j\leq i;\;m_{\!i,j\!}=0,\; i< j,\;\br\in\brac{\mbr\backslash\{0\}}^{N-1}\right\}.
\end{equation}
Similarly, we define {\bf Upper-Collinear Triangular Matrix}(U-CoLT matrix), which is a strictly upper triangular matrix:
\begin{equation}\label{eqn:ucoltdef}
    \mathcal{C}_U^N\!=\!\left\{M\in\mbr^{N\times N}\;\vert\; m_{\!i-1,j\!}/m_{\!i,j\!}=r'_{i-1},\; i< j;\;m_{\!i,j\!}=0,\; i\geq j,\;\br'\in\brac{\mbr\backslash\{0\}}^{N-2}\right\}.
\end{equation}
We call the vectors $r$ and $r'$ in \eqref{eqn:coltdef}-\eqref{eqn:ucoltdef} the \textbf{ratio vectors} of the collinear triangular matrix.
\end{definition}

The matrices introduced in ~\cref{def:colt} are termed as collinear triangular matrices (CoLT), due to the following collinearity between columns: 
\begin{equation*}
	m_{i,j}/m_{i,j+1}=m_{k,j}/m_{k,j+1}\iff m_{i,j}/m_{k,j}=m_{i,j+1}/m_{k,j+1}.
\end{equation*}

\begin{theorem}[Vector Representation of Collinear Triangular Matrix]
\label{thm:c1}
    Any L-CoLT matrix $M_L$ can be represented by its diagonal elements $\bgamma$ and the ratio vector $\br$ in Equation \eqref{eqn:coltdef}. Any U-CoLT matrix $M_U$ can be represented by its superdiagonal elements $\bgamma'$ and the ratio vector $\br'$ in \eqref{eqn:ucoltdef}.
\end{theorem}

\begin{proof}
For any L-CoLT matrix $M_L\in\mathcal{C}_L^N$, if its corresponding $\bgamma$ and $\br$ are given, then $m_{i,j}=\gamma_j\prod\limits_{k=i}^{N-1}r_k$. The proof of U-CoLT is similar.
\end{proof}

In the following, we use $\Call{L-CoLT}{\bgamma,\;\br}$, $\bgamma\in\mbr^{N},\;\br\in\mbr^{N-1}$ and $\Call{U-CoLT}{\bgamma',\;\br'}$, $\bgamma'\in\mbr^{N-1},\;\br'\in\mbr^{N-2}$ to represent a L-CoLT matrix and a U-CoLT matrix, respectively. A
specific correspondence of L-CoLT and U-CoLT is shown as follow:

For the L-CoLT matrix $M_L$
\begin{equation*}
\renewcommand\arraystretch{1.5}
	M_L = \Call{L-CoLT}{\bgamma,\;\br} = 
	\begin{pmatrix}
		\gamma_1 &   &   &   &  \\
		\gamma_1r_1& \gamma_2 &  &   &  \\
		\gamma_1r_1r_2& \gamma_2r_2&  \gamma_3  &  & \\
		\vdots &\vdots & \vdots  &\ddots &  \\
		\gamma_1\!\prod\limits_{i=1}^{\!N\!-\!1\!}\!r_i& \gamma_2\!\prod\limits_{i=2}^{\!N\!-\!1\!}\!r_i&\gamma_3\!\prod\limits_{i=3}^{\!N\!-\!1\!}\!r_i& \cdots&\;\;\gamma_{N}
	\end{pmatrix}
\end{equation*}

Similarly, for the U-CoLT matrix $M_U$
\begin{equation*}
\renewcommand\arraystretch{1.5}
	M_U = \Call{U-CoLT}{\bgamma',\;\br'} =
	\begin{pmatrix}
		0\quad &\;\gamma'_{1} &\;\gamma'_2r'_1 &\;\cdots  &\; \gamma'_{\!N\!-\!1\!}\!\prod\limits_{i=1}^{\!N\!-\!2\!}\!r'_{i} \\
		\;&\; 0 &\;\gamma'_{2} &\;\cdots  &\; \gamma'_{\!N\!-\!1\!}\!\prod\limits_{i=2}^{\!N\!-\!2\!}\!r'_{i} \\
		\;&\; &\quad\ddots &\quad\ddots  &\;\vdots \\
		\;&\; &\;  &\quad0  &\;\gamma'_{\!N\!-\!1\!}\\
		\;&\; &\;  &\; &0\; \\
	\end{pmatrix}
\end{equation*}

\medskip

The special nature of the L-CoLT and U-CoLT matrices allows us to compute matrix-vector multiplications in $O(N)$ operations.

For any $M_L=\Call{L-CoLT}{\bgamma,\;\br}$ and vector $\by\in\mbr^N$, the matrix-vector multiplication $M_L\by$ is written as
\begin{equation}
	\label{eqn:mv_summation}
	M_L\by =
	\begin{pmatrix}
		\gamma_{1}y_1 &+& 0 &+& 0 & \cdots &+& 0\\
		\gamma_1r_1y_1 &+& \gamma_2y_2 &+& 0 &\cdots  &+& 0\\
		\gamma_1r_1r_2y_1&+& \gamma_2r_2y_2&+&  \gamma_3y_3  &\cdots &+&0 \\
		\vdots &\vdots & \vdots &\vdots & \vdots & \ddots &\vdots & \vdots\\
		\gamma_1\prod\limits_{i=1}^{N-1}r_iy_1&+&\gamma_2\prod\limits_{i=2}^{N-1}r_iy_2&+& \gamma_3\prod\limits_{i=3}^{N-1}r_iy_3 & \cdots& +&\gamma_Ny_N
	\end{pmatrix}.
\end{equation}
Denote $\bp_k$ as the summation of the $k$-th row in \eqref{eqn:mv_summation}, then one has
\begin{equation*}
	\begin{aligned}
		p_1 = \gamma_{1}y_1, \quad 
		p_{k} = r_{k-1}p_{k-1} +\gamma_{k}y_{k}, \quad k=2,\cdots,N.
	\end{aligned}
\end{equation*}
Based on this recursion formula, a fast implementation is proposed in \cref{alg:LC}.
\begin{algorithm} 
    \caption{Fast L-CoLT Matrix-Vector multiplication}\label{alg:LC}
    \hspace*{0.02in} {\bf Input:} input vector $\by$ of size $N$,\;input matrix $M_L=\Call{L-CoLT}{\bgamma,\;\br}$\\
	\hspace*{0.02in} {\bf Output:} $\bp = M_L\by$
    \begin{algorithmic}[1]
        \Procedure{LCMV}{$\by,\;\bgamma,\;\br$}
        \State $p_1=\gamma_1y_1$
        \For{$i=1\;:\;N-1$}
            \State $p_{i+1} = r_{i} p_{i} +\gamma_{i+1}y_{i+1}$
        \EndFor
       
        \noindent\Return $\bp$
        \EndProcedure
    \end{algorithmic}
\end{algorithm}

Similarly, the fast matrix-vector multiplication for U-CoLT matrices is shown in \cref{alg:UC}.
\begin{algorithm} 
    \caption{Fast U-CoLT Matrix-Vector multiplication}\label{alg:UC}
    \hspace*{0.02in} {\bf Input:} input vector $\by$ of size $N$,\;input matrix $M=\Call{U-CoLT}{\bgamma',\;\br'}$\\
	\hspace*{0.02in} {\bf Output:} $\bq = M_U\by$
    \begin{algorithmic}[1]
        \Procedure{UCMV}{$\by,\;\bgamma',\;\br'$}
        \State $q_N=0$, $q_{N-1}=\gamma'_{N-1}y_{N}$
        \For{$i=2\;:\;N-1$}
            \State $q_{N-i} = r'_{N-i} q_{N-i+1} +\gamma'_{N-i}y_{N-i+1}$
        \EndFor
        
        \noindent\Return $\bq$
        \EndProcedure
    \end{algorithmic}
\end{algorithm}

Next, we denote the set $\mathcal{C}^N$ as the direct sum of $\mathcal{C}_L^N$ and $\mathcal{C}_U^N$, defined by 

\begin{definition}
\begin{equation}\label{eqn:sumset}
    \mathcal{C}^N = \mathcal{C}_L^N+\mathcal{C}_U^N=\left\{A+B\;|\;A\in\mathcal{C}_L^N,\;B\in\mathcal{C}_U^N\right\}.
\end{equation}
\end{definition}
Due to the linearity of matrix-vector multiplication, we can further develop the fast matrix-vector multiplication algorithm for matrices in $\mathcal{C}^N$, which is given in \cref{alg:recursion1d}.
\begin{algorithm} 
    \caption{CoLT Matrix-Vector multiplication}\label{alg:recursion1d}
    \hspace*{0.02in} {\bf Input:} input vector $\bx$ of size $N$, diagonal elements $\bgamma, \; \bgamma'$ and the ratio vector $\br,\;\br'$ \\
	\hspace*{0.02in} {\bf Output:} $\bp+\bq = M\by$
    \begin{algorithmic}[1]
        \Procedure{CMV}{$\by,\;\mathbf{c}^L,\;\mathbf{c}^U,\;\bgamma,\;\bgamma'$}
        \State $\bp=\Call{LCMV}{\by,\;\bgamma,\;\br}$
        \State $\bq=\Call{UCMV}{\by,\;\bgamma',\;\br'}$
        
        \noindent\Return $\bp+\bq$
        \EndProcedure
    \end{algorithmic}
\end{algorithm}
    
The space and time complexities of these algorithms are $O(N)$, which is much better than the original matrix-vector multiplication.

\subsection{Some Basic Properties}
\label{subsec:basicproperties}
In this subsection, we justify some basic properties of those matrices involved, which will be used in our algorithm.
\begin{theorem}
\label{thm:colt}
    $\brac{\mathcal{C}_L^N,\;\odot}$ and $\brac{\mathcal{C}_U^N,\;\odot}$ are Abelian groups, where $\odot$ is the Hadamard product.
    \end{theorem}
 
 
    \begin{proof} We only prove the theorem for $\brac{\mathcal{C}_L^N,\;\odot}$. It suffices to show that $\brac{\mathcal{C}_L^N,\;\odot}$ has the following properties:

Closure: For any two matrices $A=\Call{L-CoLT}{\hat{\bgamma},\;\hat{\br}}$ and $B=\Call{L-CoLT}{\tilde{\bgamma},\;\tilde{\br}}$, we set $D=A\odot B$. Since $d_{i,j}=a_{i,j}b_{i,j}$, one has
        \begin{equation}\label{eqn:closure}
    		\begin{aligned}
    			d_{i,j}/d_{i+1,j}=\brac{a_{i,j}b_{i,j}}/\brac{a_{i+1,j}b_{i+1,j}}=\hat{r}_i\tilde{r}_i,&\text{ }j=1,2,\cdots,i,
    		\end{aligned}
	    \end{equation}
	    and the strictly upper triangle part of $D$ is obviously 0, which means $D=\Call{L-CoLT}{\hat{\bgamma}\odot\tilde{\bgamma},\;\hat{r}\odot \tilde{r}} \in \mathcal{C}_L^N$.
	 
 Identity and Inverses: Let 
        \begin{equation}\label{eqn:idelement}
         E = 
             \begin{pmatrix}
                1 &   &   &   &  \\
        		1& \quad1 &  &   &  \\
        		1& \quad1&  \quad1  &  & \\
        		\vdots &\quad\vdots & \quad\vdots  &\quad\ddots &  \\
        		1& \quad1& \quad1& \quad\cdots&\quad1
             \end{pmatrix}\in\mathcal{C}_L^N,
         \end{equation}
         then for any $A=\Call{L-CoLT}{\bgamma,\;\br}$, $A\odot E=E\odot A=A$, which means $E$ is the identity element. Let \begin{equation*}
         B = 
             \begin{pmatrix}
                1/a_{11} &   &   &   &  \\
        		1/a_{21}& \quad 1/a_{22} &  &   &  \\
        		1/a_{31}& \quad1/a_{32}&  \quad1/a_{33}  &  & \\
        		\vdots &\quad\vdots & \quad\vdots  &\quad\ddots &  \\
        		1/a_{n1}& \quad1/a_{n2}& \quad1/a_{n3}& \quad\cdots&\quad1/a_{nn}
             \end{pmatrix}.
         \end{equation*}Since \[b_{i,j}/b_{i+1,j}=a_{i+1,j}/a_{i,j}=1/r_i,\;j=1,2,\cdots,i, \]then $B\in\mathcal{C}_L^N$. And obviously, $A\odot B=B\odot A=E$, which means $B$ is the inverse of $A$.
	 
Commutativity and Associativity:  The commutativity and associativity can be derived from the commutative and associative law of real number multiplication.

\end{proof}
\medskip
Based on the above theorem, we can deduce directly 
\begin{corollary}\label{cor:colt1}
     $\brac{\mathcal{C}^N,\;\odot}$ is an abelian group with identity element $\mathbf{1}_N\mathbf{1}_N^T$.
\end{corollary}
\medskip
\begin{theorem}
\label{thm:colt2}
    For any vector $\bx\in\brac{\mbr\backslash\{0\}}^{N}$, $f_{\bx}\brac{M}=\brac{\text{diag}\brac{\bx}} M$ and $g_{\bx}\brac{M}=M\brac{\text{diag}\brac{\bx}}$ are permutations in $\mathcal{C}_L^N$ and $\mathcal{C}_U^N$.
\end{theorem}
\begin{proof} We only prove the theorem for $\mathcal{C}_L^N$.
    
Closure: For any vector $\bx\in\brac{\mbr\backslash\{0\}}^{N}$, and $M_L=\Call{L-CoLT}{\bgamma,\;\br}$, let $E$ be the one defined in Equation \eqref{eqn:idelement}. Since
\[       E_1 = \brac{\text{diag}\brac{\bx}}E=
             \begin{pmatrix}
                x_{1} &   &   &   &  \\
        		x_{2}& \; x_{2} &  &   &  \\
        		x_{3}& \; x_{3}&  \; x_{3}  &  & \\
        		\vdots &\;\vdots & \;\vdots  &\;\ddots &  \\
        		x_{n}& \; x_{n}& \; x_{n}& \;\cdots&\; x_{n}
             \end{pmatrix}\in\mathcal{C}_L^N,\]
\[        E_2 = E\brac{\text{diag}\brac{\bx}}=
         \begin{pmatrix}
                x_{1} &   &   &   &  \\
        		x_{1}& \; x_{2} &  &   &  \\
        		x_{1}& \; x_{2}&  \; x_{3}  &  & \\
        		\vdots &\;\vdots & \;\vdots  &\;\ddots &  \\
        		x_{1}& \; x_{2}& \; x_{3}& \;\cdots&\; x_{n}
            \end{pmatrix}\in\mathcal{C}_L^N,\]
        we have 
        \begin{equation}\label{eqn:lemma_right}
    		\begin{aligned}
    			f_{\bx}\brac{M}=\brac{\text{diag}\brac{\bx}}M=\brac{\text{diag}\brac{\bx}}E\odot M = E_1\odot M\in\mathcal{C}_L^N;\\
    			g_{\bx}\brac{M}=M\brac{\text{diag}\brac{\bx}}=M\odot E\brac{\text{diag}\brac{\bx}} = M\odot E_2\in\mathcal{C}_L^N.
    		\end{aligned}
	    \end{equation}
	    The last set membership can be derived by the closure of $\brac{\mathcal{C}_L^N,\;\odot}$ proved in \cref{thm:colt}. Hence, $f_{\bx}$ and $g_{\bx}$ are maps from $\mathcal{C}_L^N$ to itself.
	    
Injectiveness: For any two matrices $A=\Call{L-CoLT}{\hat{\bgamma},\;\hat{r}}$ and $B=\Call{L-CoLT}{\tilde{\bgamma},\;\tilde{\br}}$, let $D_1$ be the inverse of $E_1$ and $D_2$ be the inverse of $E_2$. If $f_{\bx}\brac{A}=f_{\bx}\brac{B}$, then \[A = D_1\odot E_1\odot A=D_1\odot f_{\bx}\brac{A} =D_1\odot f_{\bx}\brac{B}=D_1\odot E_1\odot B = B.\]If $g_{\bx}\brac{A}=g_{\bx}\brac{B}$, then \[A = A\odot E_2\odot D_2= g_{\bx}\brac{A}\odot D_2 =g_{\bx}\brac{B}\odot D_2=B\odot E_2\odot D_2 = B,\] which means $f_{\bx}\brac{\cdot}$ and $g_{\bx}\brac{\cdot}$ are injective functions.
	    
Surjectiveness: For any $M\in\mathcal{C}_L^N$, let $Q_1=D_1\odot M$ and $Q_2=M\odot D_2$, then
\begin{equation*}
\begin{aligned}
    f_{\bx}\brac{Q_1}=E_1\odot D_1\odot M=E\odot M=M; \\
    g_{\bx}\brac{Q_2}=M\odot D_2\odot E_2 = M\odot E = M,
\end{aligned}
\end{equation*} which means $f_{\bx}\brac{\cdot}$ and $g_{\bx}\brac{\cdot}$ are surjective functions.
\end{proof}
\begin{corollary}\label{cor:colt2}
    $f_{\bx}\brac{\cdot}$ and $g_{\bx}\brac{\cdot}$ are permutations in $\mathcal{C}^N$.
\end{corollary}

\section{The Fast Sinkhorn II} \label{sec:fs2}

In this section, we will discuss the implementation details to accelerate IPOT. In \cref{alg:IPOT}, three parts lead to $O(N^2)$ algorithm complexity, i.e., the matrix Hadamard product (line 4), the matrix-vector multiplication (lines 6-7), and the matrix scaling (line 8). They all rely on the representation and manipulation of L/U CoLT matrices.

    For two discrete distributions on a 1D uniform mesh grid with a grid spacing of $h$, by introducing the notation $\lambda = e^{-h/\delta}$, the kernel matrix $K$ is written as
    \begin{equation}\label{eqn:originalK}
		K =
		\begin{pmatrix}
			1 & \lambda & \lambda^2 & \cdots &\lambda^{N-1}\\
			\lambda& 1 &\lambda &\cdots  &\lambda^{N-2}\\
			\lambda^2& \lambda&  1  &\cdots &\lambda^{N-3}\\
			\vdots &\vdots & \vdots &\ddots & \vdots\\
			\lambda^{N-1}& \lambda^{N-2}& \lambda^{N-3}& \cdots&1
		\end{pmatrix} \in \mathcal{C}^N.
	\end{equation}
	Below we discuss step by step of IPOT (\cref{alg:IPOT}) to reduce the complexity:
	\begin{itemize}
		\item line 2: $\Gamma=\mathbf{1}_N\mathbf{1}^T_N\in \mathcal{C}^N$, we only need four vectors $(\bgamma,\bgamma',\br,\br')$ to represent $\Gamma$ according to \cref{thm:c1}. 
		\item line 4: the matrix Hadamard product $Q=K\odot\Gamma \in \mathcal{C}^N$ according to \cref{thm:colt} and \cref{cor:colt1}. By updating the four representation vectors $(\bgamma,\bgamma',\br,\br')$, we can obtain $Q$ with $O(N)$ cost.
		\item lines 6-7: the matrix-vector multiplication $Q^T\bphi$ and $Q\bpsi$ can be implemented with $O(N)$ cost according to \cref{alg:recursion1d}.
		\item line 8: the matrix scaling $\Gamma=\text{diag}(\bphi)Q\text{diag}(\bpsi)\in \mathcal{C}^N$ according to \cref{thm:colt2} and \cref{cor:colt2}. By updating the four representation vectors $(\bgamma,\bgamma',\br,\br')$, we can obtain $\Gamma$ with $O(N)$ cost.
	\end{itemize}
	
	Based on the above discussions, we proposed the FS-2 algorithm with $O(N)$ complexity. The pseudo-code is presented in \cref{alg:FS-2}.
    \begin{algorithm} 
		\caption{1D FS-2 Algorithm}
		\label{alg:FS-2}
		\hspace*{0.02in} {\bf Input:} $\bu,\;\bv\in\mbr^{N}$;\ $L,\mathrm{itr\_max}\in\mbn^+$;\ $h,\delta\in\mbr$\\
		\hspace*{0.02in} {\bf Output:} $W_1(\bu,\bv)$
		\begin{algorithmic}[1]
			\State $\lambda\gets e^{-h/\delta}$;\  $\bphi,\bpsi\gets\frac{1}{N}\mathbf{1}_N$;\  $\br,\bs\gets\mathbf{0}_{N}$
			\State $\balpha^{L},\;\bbeta^{L},\;\balpha^{U},\;\bbeta^{U}\gets\lambda\mathbf{1}_{N-1}$;\;$\bgamma \gets \mathbf{1}_N$;\;$\bgamma' \gets \lambda\mathbf{1}_{N-1}$
			\For{$t=1\;:\;\mathrm{itr\_max}$}
			\For{$\ell=1\;:\;L$}
			\State $\br\gets \Call{CMV}{\bphi,\bbeta^{L},\bbeta^{U},\bgamma,\bgamma'}$
			\State $\bpsi\gets\bv\oslash\br$
			\State $\bs\gets \Call{CMV}{\bpsi,\balpha^{L},\balpha^{U},\bgamma,\bgamma'}$
			\State $\bphi\gets\bu\oslash\bs$
			\EndFor
			\For{$i=1\;:\;N-1$}
			\State $\alpha^L_i\gets\lambda\alpha^L_i\brac{\phi_{i+1}/\phi_i}$,\;$\beta^L_i\gets\lambda\beta^L_i\brac{\psi_{i+1}/\psi_i}$
			\State $\gamma'_{i}\gets\lambda\gamma'_{i}\phi_{i}\psi_{i+1}$
			\EndFor
			\For{$i=1\;:\;N-2$}
			\State $\alpha^U_i\gets\lambda\alpha^U_i\brac{\phi_{i}/\phi_{i+1}}$,\;$\beta^U_i\gets\lambda\beta^U_i\brac{\psi_{i}/\psi_{i+1}}$
			\EndFor
			\State $\bgamma\gets\bphi\odot\bpsi\odot\bgamma$
			\EndFor		
			\noindent\Return $W_1(\bu,\bv)$
		\end{algorithmic}
	\end{algorithm}	
	
	There is a minor flaw in the above algorithm. The computational cost of $W_1(\bu,\bv)$ is still $O(N^2)$ in the last step. This was also ignored in our previous paper~\cite{liao2022fast}. Now, we would like to discuss this issue. The computation of $W_1\brac{\bu,\bv}=\left\langle C,\;\Gamma\right\rangle$ can be regarded as the summation of all elements of the following matrix.

    \begin{equation}
    	\label{eqn:w1_summation}
    	\renewcommand\arraystretch{1.3}
    	\small
    	C\odot\Gamma=
    	\begin{pmatrix}
    		0&h\gamma'_{1} &2h\gamma'_2r'_1 &\cdots  & (\!N\!-\!1\!)h\gamma'_{\!N\!-\!1\!}\!\prod\limits_{i=1}^{\!N\!-\!2\!}\!r'_{i}\\
    		h\gamma_1r_1 &0 &h\gamma'_{2} &\cdots  &(\!N\!-\!2\!)h\gamma'_{\!N\!-\!1\!}\!\prod\limits_{i=2}^{\!N\!-\!2\!}\!r'_{i} \\
    		2h\gamma_1r_1r_2 &h\gamma_2r_2 &0 &\cdots &(\!N\!-\!3\!)h\gamma'_{\!N\!-\!1\!}\!\prod\limits_{i=3}^{\!N\!-\!2\!}\!r'_{i}\\
    		\vdots & \vdots & \vdots & \ddots & \vdots\\
    		\brac{\!N\!-\!1\!}h\gamma_1\!\prod\limits_{i=1}^{\!N\!-\!1\!}\!r_i &\;\;\brac{\!N\!-\!2\!}h\gamma_2\!\prod\limits_{i=2}^{\!N\!-\!1\!}\!r_i & \;\;\brac{\!N\!-\!3\!}h\gamma_3\!\prod\limits_{i=3}^{N\!-\!1}\!r_i & \cdots &0
    	\end{pmatrix}.
    \end{equation}
        We separate the summation of the matrix to the lower and strictly upper triangular parts. Thus, the $k$-th line summation of two parts can be written as
        \begin{equation*}
        		p_k=\sum\limits_{i=1}^k \omega_{ki},\quad
		q_k=\sum\limits_{i=k+1}^{N} \omega_{ki}.
        \end{equation*}
        We can consider the following recursive computation
    	\begin{equation}\label{eqn:CGAMMA}
        		\begin{aligned}
            		& p_1 = 0,\quad p_2 = h\gamma_1r_1,\quad p'_2 = h\gamma_1r_1 + h\gamma_2,\\
            		& p_i = r_{i-1}\brac{p_{i-1} + p'_{i-1}},\quad p'_{i} = r_{i-1}p'_{i-1}+h\gamma_{i},\quad i = 3,\;4,\cdots,\;N. \\
            		& q_{N} = 0, \quad q_{N-1}=h\gamma'_{N-1},\quad q'_{N-1}=h\gamma'_{N-1} + h, \\
            		& q_{j}=r'_{j}\brac{q_{j+1}+q'_{j+1}},\quad q'_{j}=r'_{j}q'_{j+1}+h,\quad j = 1,\;2,\cdots,\;N-2.
		\end{aligned}
	\end{equation}
	Thus, the Wasserstein-1 metric can be finally obtained with $O(N)$ cost
	\begin{equation*}
		W_1(\bu,\bv)=\left\langle C,\;\Gamma\right\rangle=\sum_{i=1}^N\brac{p_i+q_i}.
	\end{equation*} 
        
    
\section{Extension to high dimension}
\label{sec:fs2hd}
In this section, we illustrate how the FS-2 algorithm generalizes to higher dimensions using the two-dimensional case as an example. 
\subsection{Block Collinear Triangular Matrix}
Hereinafter, for $A\in\mbr^{MN\times MN}$, we break it into $M^2$ uniform blocks with size $N\times N$:
    \begin{equation*}
    \renewcommand\arraystretch{2}
        	A = \left(\begin{array}{c|c|c|c|c}
    		A_{1,1} & A_{1,2} & A_{1,3} & \cdots & A_{1,M} \\
    		\hline
    		A_{2,1} & A_{2,2} & A_{2,3} & \cdots & A_{2,M} \\
    		\hline
    		\vdots &\vdots &\vdots  &\ddots & \vdots \\
    		\hline
    		A_{M,1} & A_{M,2} & A_{M,3} & \cdots & A_{M,M}
    	\end{array}\right).
    \end{equation*}
And for vectors $\bx\in\mbr^{kN}$, we break it into $k$ uniform blocks as $\brac{\bx_1,\bx_2,\cdots,\bx_k}^T$, in which
\begin{equation*}
	\bx_i=\brac{x_{1+(i-1)N},\;x_{2+(i-1)N},\;\cdots,\;x_{iN}}^T,\quad i=1,2,\cdots,k.
\end{equation*}To carry out the fast implementation of the matrix-vector multiplication of the block matrix above, we generalize the definition of $\mathcal{C}^{N}$ in \eqref{eqn:sumset} to the block case as
\medskip
\begin{definition}
    \begin{equation*}
    \begin{aligned}
        \mathcal{C}^{N,M}=\{A\in&\mbr^{MN\times MN}\vert A_{k,k}\in\mathcal{C}^N; \; \br^L,\br^U\in\brac{\mbr\backslash\{0\}}^{\brac{M-1}N};
        	\\ &A_{i+1,j}=\brac{\text{diag}\brac{\br^L_i}}A_{i,j},\; j\leq i;\;A_{i-1,j}=\brac{\text{diag}\brac{\br^U_{i-1}}}A_{i,j},\;i\leq j\}
    \end{aligned}.
    \end{equation*}
\end{definition}

Since $\mathcal{C}^{N,M}$ is a generalization of $\mathcal{C}^N$, we can also use the strategy of \cref{alg:recursion1d} in blocks to reduce the computational cost of matrix-vector multiplications. For a vector $\bx\in\mbr^{NM}$, the matrix-vector multiplication $A\bx$ is written as
\begin{equation}\label{eqn:NBS_summation 2d}
	A \bx =
	\begin{pmatrix}
		A_{1,1}\bx_1 &+& A_{1,2}\bx_2 &+& A_{1,3}\bx_3 & \cdots &+& A_{1,M}\bx_M\\
		A_{2,1}\bx_1 &+& A_{2,2}\bx_2 &+& A_{2,3}\bx_3 & \cdots &+& A_{2,M}\bx_M\\
		A_{3,1}\bx_1 &+& A_{3,2}\bx_2 &+& A_{3,3}\bx_3 & \cdots &+& A_{3,M}\bx_M \\
		\vdots &\vdots & \vdots &\vdots & \vdots & \ddots &\vdots & \vdots\\
		A_{M,1}\bx_1 &+& A_{M,2}\bx_2 &+& A_{M,3}\bx_3 & \cdots &+& A_{M,M}\bx_M
	\end{pmatrix}.
\end{equation}
We separate the summation of row $k$ to the lower triangular part $\bp_k$ and the strictly upper triangular part $\bq_k$. Then computing $A \bx$ is formulated as
\begin{equation*}\label{meshpoly 2d}
	A \bx = \bp+\bq, \quad
	\bp_k = \sum_{i=1}^{k} A_{k,i}\bx_i, \quad
	\bq_k = \sum_{i=k+1}^{M} A_{k,i}\bx_i, \quad k=1,\cdots,M.
\end{equation*}
If $A$ is in $\mathcal{C}^{N,M}$ with $\br^L$ and $\br^U$, instead of directly calculating $\bp_k$ and $\bq_k$, a successive computation is used
\begin{equation}\label{eqn:recursion2d}
	\begin{aligned}
		&\bp_1 = A_{1,1}\bx_1, \quad 
		\bp_{k} = \br^L_{k-1}\odot\bp_{k-1} + A_{k,k}\bx_{k}, \quad k=2,\cdots,M, \\
		&\bq_M = \mathbf{0}_N, \quad
		\bq_{k-1} = \br^U_{k-1}\odot(\bq_{k}+A_{k,k}\bx_{k}), \quad k=M,M-1,\cdots,2.
	\end{aligned}
\end{equation}
Since the computation of $A_{k,k}\bx_k$ can be carried out with $O(N)$ complexity by using \cref{alg:recursion1d}, the whole computation is of $O(NM)$ complexity.

Similar to \cref{thm:colt} and \cref{thm:colt2}, $\mathcal{C}^{N,M}$ is closed under the Hadamard product and matrix scaling.
\begin{theorem}\label{thm:2d} $\brac{\mathcal{C}^{N,M},\odot}$ is an Abelian group; Matrix scaling operations are permutations in $\mathcal{C}^{N,M}$.
    \begin{proof}
        For any $A\in\mathcal{C}^{N,M}$, since the diagonal blocks of $A$ are in $\mathcal{C}^N$, by \cref{cor:colt2}, all blocks in $A$ are in $\mathcal{C}^N$. Then the two properties can be proved in a similar way as in \cref{subsec:basicproperties}.
    \end{proof}
\end{theorem}

\subsection{The 2D FS-2 Algorithm}
Consider two discretized probabilistic distributions 
\begin{align*}
    \bu &= \brac{u_{11},u_{21},\cdots,u_{N1},u_{12},\cdots,u_{i_1j_1}, \cdots, u_{NM}},\\
    \bv &= \brac{v_{11},v_{21},\cdots,v_{N1},v_{12},\cdots, v_{i_2j_2}, \cdots, v_{NM}},
\end{align*}
on a uniform 2D mesh of size $N\times M$ with a vertical spacing of $h_1$ and a horizontal spacing of $h_2$. The corresponding kernel matrix is written as
\begin{equation*}
\renewcommand\arraystretch{2}
    	K = \left(\begin{array}{c|c|c|c|c}
		K_0 & \lambda_2 K_0 & \lambda_2^2K_0 & \cdots & \lambda_2^{M-1}K_0 \\
		\hline
		\lambda_2 K_0& K_0 & \lambda_2 K_0 &\cdots & \lambda_2^{M-2}K_0 \\
		\hline
		\vdots &\vdots &\vdots  &\ddots & \vdots \\
		\hline 
		\lambda_2^{M-1}K_0\; & \lambda_2^{M-2}K_0\; & \lambda_2^{M-3}K_0\; & \cdots & K_0
	\end{array}\right),
\end{equation*}
where the sub-matrix
\begin{equation*}
    K_0=\begin{pmatrix}
	1 & \lambda_1 &\cdots &\lambda_1^{N-1} \\
	\lambda_1 & 1 &\cdots &\lambda_1^{N-2} \\
	\vdots	& \vdots &\ddots &\vdots	\\
	\lambda_1^{N-1} & \lambda_1^{N-2} & \cdots & 1
\end{pmatrix},
\end{equation*}
and
\begin{equation*}
	\lambda_1=e^{-h_1/\delta}, \quad \lambda_2 = e^{-h_2/\delta}.
\end{equation*}

Obviously, the original 2D kernel $K$ contains blocks which are multiples of the 1D kernel, hence belongs to $\mathcal{C}^{N,M}$. By an analysis similar to that in \cref{sec:fs2} and using \cref{thm:2d}, the matrices $Q$ and $\Gamma$ are in $\mathcal{C}^{N,M}$ throughout the course of the iteration, which means that all the matrix-vector multiplications can be carried out by using recursion \eqref{eqn:recursion2d}. Thus, the total cost of matrix-vector multiplication of our FS-2 algorithm for 2D Wasserstein-1 metric is reduced to $O(NM)$. 

In the 2D FS-2 algorithm, we use ` $\widehat{}$ ’ to distinguish whether it is a coefficient of the block or the inner sub-matrix, and update them simultaneously after an inner loop. The pseudo-code is presented in \cref{alg:FS-22d}. Updating the coefficients of inner sub-matrices and computation of $W_1\brac{\bu,\bv}$ are omitted in the pseudo-code since they are similar to the 1D case, which we have described in detail.
	\begin{algorithm} 
		\caption{2D FS-2 Algorithm}
		\label{alg:FS-22d}
		\hspace*{0.02in} {\bf Input:} $\bu,\;\bv\in\mbr^{NM}$;\ $L,\mathrm{itr\_max}\in\mbn^+$;\ $h_1,h_2,\delta\in\mbr$\\
		\hspace*{0.02in} {\bf Output:} $W_1(\bu,\bv)$
		\begin{algorithmic}[1]
			\State $\lambda_1\gets e^{-h_1/\delta}$;\;$\lambda_2\gets e^{-h_2/\delta}$;\; $\bphi^{(0)},\bpsi^{(0)} \gets \frac{1}{NM}\mathbf{1}_{NM}$; $\bp^{(0)},\br^{(0)},\bq^{(0)},\bs^{(0)}\gets\mathbf{0}_{NM}$
			\State $\bgamma \gets \mathbf{1}_{NM},\;\bgamma',\;\balpha^{L},\;\bbeta^{L},\lambda_1\mathbf{1}_{(N-1)M},\;\balpha^{U},\;\bbeta^{U}\gets\lambda_1\mathbf{1}_{(N-2)M}$
			\State $\widehat{\balpha}^{L},\;\widehat{\bbeta}^{L},\;\widehat{\balpha}^{U},\;\widehat{\bbeta}^{U}\gets\lambda\mathbf{1}_{N(M-1)}$
			\While{$t=1\;:\;\mathrm{itr\_max}$}
			\For{$\ell=1\;:\;L$}
			\State $\br_{1} \gets\Call{CMV}{\bphi_1,\;\bbeta_1}$,\;$s_{M} \gets \mathbf{0}_N$			
			\For{$i=1\;:\;M-1$}
			\State $\br_{i+1} \gets \widehat{\bbeta}^L_{i}\odot \br_{i}+\Call{CMV}{\bphi_{i+1},\bbeta_{i+1}^L,\bbeta_{i+1}^U,\bgamma_{i+1},\bgamma'_{i+1}}$
			\State $\bs_{\!N-i\!} \!\gets\! \widehat{\bbeta}^U_{\!N-i\!}\!\odot\!(\bs_{\!N-i+1\!}+\Call{\!CMV\!}{\bphi_{\!N-i+1\!},\bbeta_{\!N-i+1\!}^L,\bbeta_{\!N-i+1\!}^U,\bgamma_{\!N-i+1\!},\bgamma'_{\!N-i+1\!}})$
			\EndFor
			\State $\bpsi\gets\bv\oslash\brac{\br+\bs}$
			\State $\bp_{1} \gets\Call{CMV}{\bpsi_1,\;\balpha_1}$, $\bq_{M} \gets \mathbf{0}_N$
			\For{$i=1\;:\;M-1$}
			\State $\bp_{i+1} \gets \widehat{\balpha}^L_{i}\odot \bp_{i}+\Call{CMV}{\bpsi_{i+1},\balpha_{i+1}^L,\balpha_{i+1}^U,\bgamma_{i+1},\bgamma'_{i+1}}$
			\State $\bq_{\!N-i\!}\!\gets\! \widehat{\balpha}^U_{\!N-i\!}\!\odot\!(\bq_{\!N-i+1\!}+\Call{\!CMV\!}{\bpsi_{\!N\!-\!i\!+\!1\!},\balpha^L_{\!N-i+1\!},\balpha^U_{\!N-i+1\!},\bgamma_{\!N-i+1\!},\bgamma'_{\!N-i+1\!}})$
			\EndFor
			\State $\bphi\gets\bu\oslash\brac{\bp+\bq}$
			\EndFor
			\For{$i=1\;:\;M-1$}
			\State $\widehat{\balpha}^L_{i} \gets \lambda_2\brac{\bphi_{i+1}/\bphi_{i}}\odot\widehat{\balpha}^L_{i}$,\;$\widehat{\bbeta}^L_{i} \gets \lambda_2\brac{\bpsi_{i+1}/\bpsi_{i}}\odot\widehat{\bbeta}^L_{i}$
			\State $\widehat{\balpha}^U_{i} \gets \lambda_2\brac{\bphi_{i}/\bphi_{i+1}}\odot\widehat{\balpha}^U_{i}$,\;$\widehat{\bbeta}^U_{i} \gets \lambda_2\brac{\bpsi_{i}/\bpsi_{i\!+\!1}}\odot\widehat{\bbeta}^U_{i}$
			\EndFor
			\State Update $\bgamma,\;\bgamma',\;\balpha^L,\;\bbeta^L,\;\balpha^U,\;\bbeta^U$
			\EndWhile{}			
			\noindent\Return $W_1(\bu,\bv)$
		\end{algorithmic}
	\end{algorithm}

\section{Numerical Experiments}
\label{sec:5}
In this section, we carry out three numerical experiments to evaluate the FS-2 algorithm, including one 1D example and two 2D examples. The true Wasserstein metric $W_{LP}$ is obtained by solving the original OT \eqref{eqn:ot} using interior-point methods~\cite{dikin1967iterative,karmarkar1984new}. In our experiments, for both IPOT and FS-2, the number of inner loops is set as $L=20$ and the regularization parameter $\delta^{(t)}$ is set to $1$. The number of iterations here is the total number of loops: $\#\mathrm{iteration}=\mathrm{itr\_max}\times L$. In order to deal with the difficulties caused by zeros, we utilize the rescaling method in~\cite{lzy2021}:
\begin{equation} \label{eqn:shift}
        	D(f,g)=W_{1}\brac{\frac{\frac{f}{\norm{f}}+\eta}{1+N\eta}, \frac{\frac{g}{\norm{g}}+\eta}{1+N\eta}},
    \end{equation}
In the following, we refer to formula \eqref{eqn:shift} for numerical stability with $\eta=10^{-5}$. All the experiments are conducted on a platform with 128G RAM, and one Intel(R) Xeon(R) Gold 5117 CPU @2.00GHz with 14 cores. 

\subsection{1D Gaussian distributions}
\label{subsec:1dgaussian}
We consider the Wasserstein-1 metric between two mixtures of 1D Gaussian distributions: $0.4\mathcal{N}\brac{60,64}+0.6\mathcal{N}\brac{40,36}$ and $0.5\mathcal{N}\brac{35,81}+0.5\mathcal{N}\brac{70,81}$, which is the experiment setting in~\cite{xie2020fast}. Input vectors $\bu$ and $\bv$ are generated by integration on the uniform discretization of interval $[0,100]$ with node size $N$. 

\begin{figure} 
    \centering
    \includegraphics[width=0.6\linewidth]{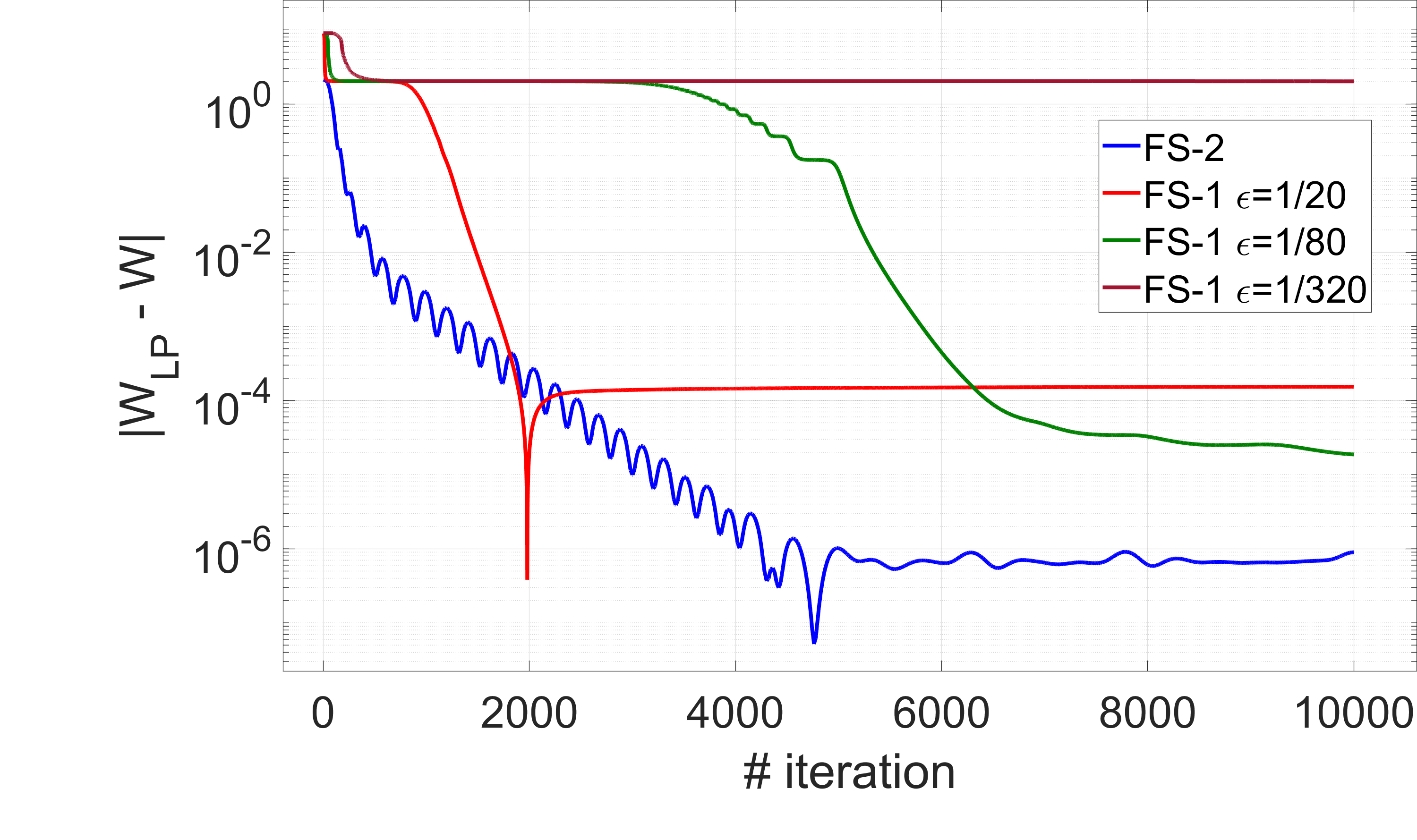}
    \caption{The 1D Gaussian distribution problem. The errors between the numerical results generated by FS-1 or FS-2 and the true Wasserstein-1 metric w.r.t. number of iterations.}
    \label{fig:1diter}
\end{figure}

We first compare the convergence of FS-1 and FS-2 for $N=1000$. We tested $100$ experiments, and each experiment was performed for $10,000$ iterations. In \cref{fig:1diter}, the differences of the Wasserstein-1 metric between the true solution $W_{LP}$ and the numerical solutions generated by FS-1 and FS-2 are depicted. As expected, as $\veps$ decreases, the error of FS-1 decreases gradually after the iterations converge. We can observe this for $\veps=1/20$ and $\veps=1/80$. However, for $\veps=1/320$, we can not observe convergence. In fact, the error does not drop over $10,000$ iterations. This is because $\veps$ is too small, making updates extremely slow. In fact, after $20,000$ iterations, the result of $\veps=1/320$ will continue to drop, and the final error is smaller than that of $\veps=1/20$ and $\veps=1/80$. However, in any case, the results of FS-1 are far inferior to those of FS-2, both in terms of accuracy and convergence rate.

\begin{figure} 
    \centering
    \includegraphics[width=\linewidth]{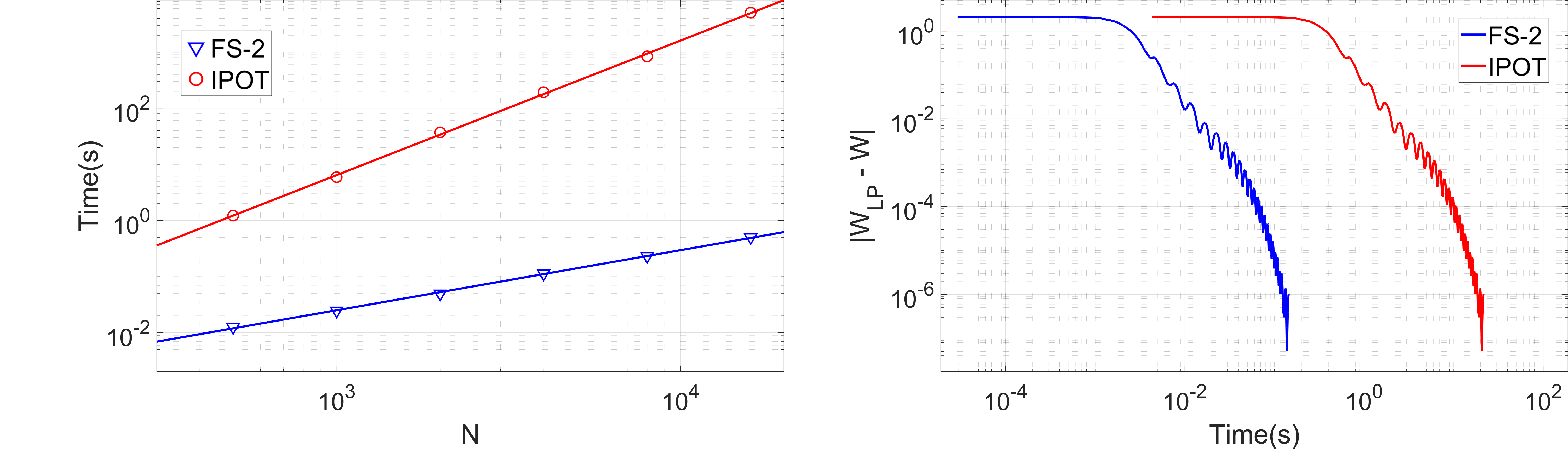}
    \caption{The 1D Gaussian distribution problem. Left: The comparison of computational time between the FS-2 algorithm and the IPOT method with different numbers of grid points $N$. Right: The computational time required to reach the absolute error of the Wasserstein-1 metric.}
    \label{fig:1dGaussian}
\end{figure}

The averaged computational time of the IPOT method and the FS-2 algorithm is given in \cref{table:1dGaussian} and \cref{fig:1dGaussian} (left). Apparently, the FS-2 algorithm has achieved an overwhelming advantage over the IPOT method in terms of computational speed, and ensures that the transport plans of the two are almost the same. This replicates the advantages of the FS-1 algorithm over the Sinkhorn algorithm. According to the data fitting results, the empirical complexity of the FS-2 algorithm is $O(N^{1.07})$, which is much smaller than the $O(N^{2.40})$ complexity of the IPOT method. At last, we show the computational time required to reach the absolute error of the Wasserstein-1 metric for $N=1,000$ in \cref{fig:1dGaussian} (right). Clearly, the FS-2 algorithm has an advantage of two orders of magnitude in computational time compared to the IPOT method.

\begin{table}
	\caption{The 1D Gaussian distribution problem. The comparison between the IPOT method and the FS-2 algorithm with the different number of grid points $N$. Columns 2-4 are the averaged computational time of the two algorithms and the speed-up ratio of the FS-2 algorithm. Column 5 is the Frobenius norm of the difference between the transport plan computed by the two algorithms.}
	\centering
	\begin{tabular}{cclcc}
		\toprule
		\multirow{2}{*} N &
		\multicolumn{2}{c}{Computational time (s)} & \multirow{2}{*}{Speed-up ratio} & \multirow{2}{*}{$\norm{P_{FS}-P}_F$} \\
		
		\cline{2-3} &FS-2 &IPOT \\
		\midrule
		500   & $1.25\times10^{-2}$ &$1.22\times10^{0}$ &$9.76\times10^{1}$   & $2.09\times10^{-15}$  \\ 
		2000  & $4.95\times10^{-2}$ & $3.73\times10^{1}$ &$7.52\times10^{2}$   & $6.65\times10^{-16}$ \\ 
		8000  & $2.32\times10^{-1}$ & $8.41\times10^{2}$ & $3.63\times10^{3}$ &
		$8.57\times10^{-16}$\\ 
		\bottomrule
	\end{tabular}
	\label{table:1dGaussian}
\end{table}

\subsection{2D Random distributions}
\label{subsec:2drandom}
\begin{figure} 
    \centering
    \includegraphics[width=0.6\linewidth]{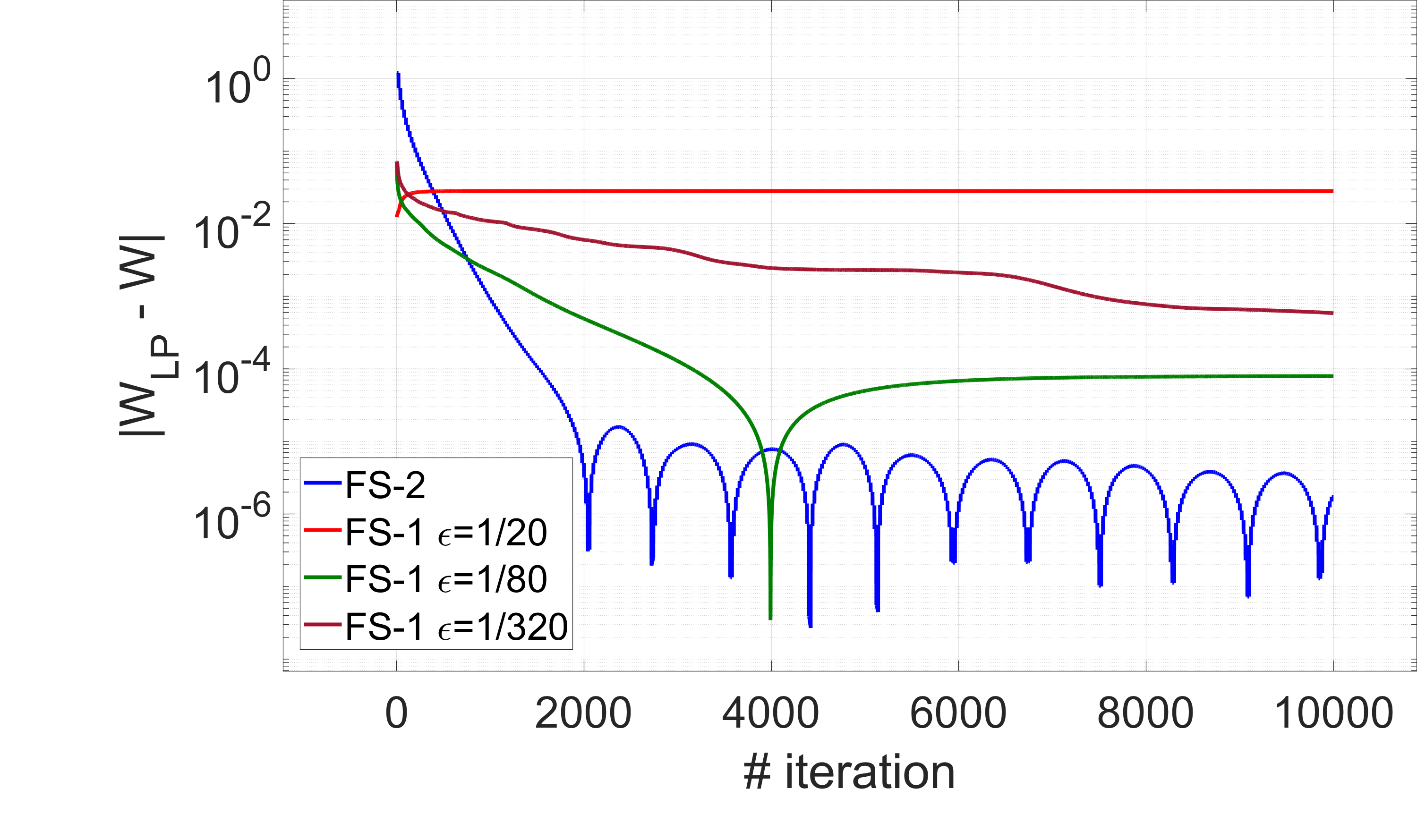}
    \caption{The 2D random distribution problem. The errors between the numerical results generated by FS-1 or FS-2 and the true Wasserstein-1 metric w.r.t. number of iterations.}
    \label{fig:2diter}
\end{figure}

Next, we compute the Wasserstein-1 metric between two $N\times N$ dimensional random vectors whose elements obey the uniform distribution on $(0,1)$. Without loss of generality, we set $h_1=h_2=0.1$. We also tested $100$ experiments, and each experiment was performed for $10,000$ iterations. We hope to test the performance of the FS-2 algorithm in 2D through this example. The differences in the Wasserstein-1 metric between the true solution $W_{LP}$ and the numerical solutions generated by FS-1 and FS-2 are shown in \cref{fig:2diter}. From this, we can observe that FS-1 converges quickly for $\veps=1/20$, but the error is large. When $\veps=1/80$, the iteration converges at about $5,000$ steps. The error keep decreasing even after $10,000$ steps for $\veps=1/320$. However, their errors and convergence speed are not as good as FS-2.

\begin{figure} 
    \centering
    \includegraphics[width=\linewidth]{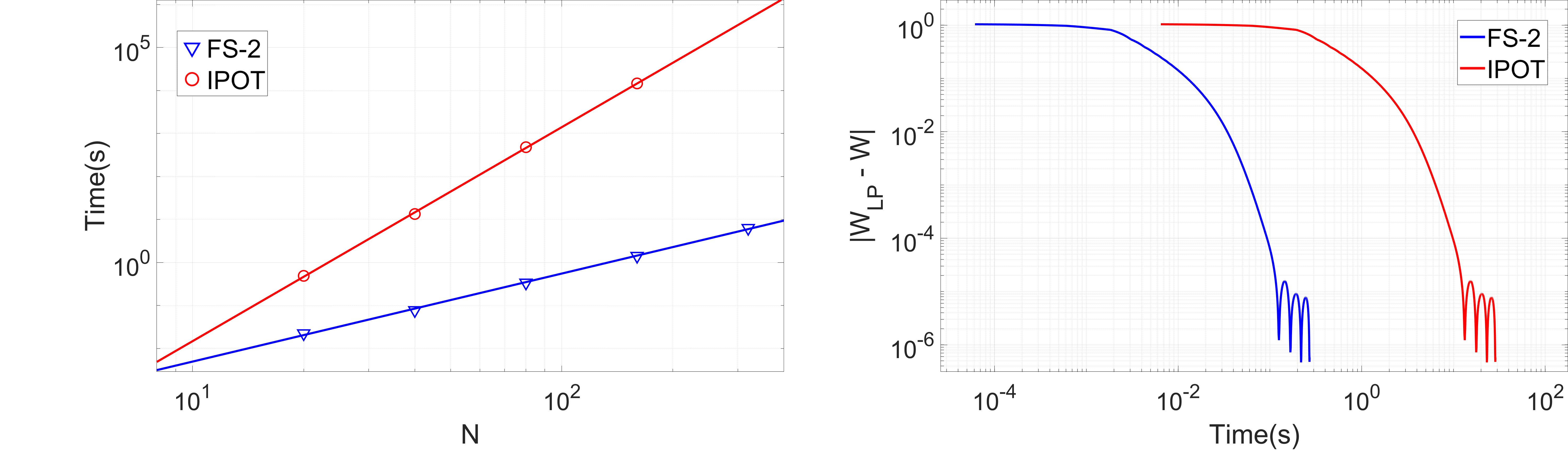}
    \caption{The 2D random distribution problem. Left: The comparison of computational time between the FS-2 algorithm and the IPOT method with different numbers of grid points $N$. Right: The computational time required to reach the absolute error of the Wasserstein-1 metric.}
    \label{fig:2drandom}
\end{figure}

The averaged computational time of the IPOT method and the FS-2 algorithm is given in \cref{table:2drandom} and \cref{fig:2drandom} (left). According to the data fitting results, the empirical complexity of the FS-2 algorithm is $O(N^{2.05})$, which is much smaller than the $O(N^{4.98})$ complexity of the IPOT method. The computational time required to reach the absolute error of the Wasserstein-1 metric for $N\times N=32\times 32$ is also presented in \cref{fig:2drandom} (right). Similar to the previous subsection, we can also observe the huge computational efficiency of the FS-2 algorithm over the IPOT method.

\begin{table}
	\caption{The 2D random distribution problem. The comparison between the IPOT method and the FS-2 algorithm with  different total number of grid nodes $N\times N$. Columns 2-4 are the averaged computational time of the two algorithms and the speed-up ratio of the FS-2 algorithm. Column 5 is the Frobenius norm of the difference between the transport plan computed by the two algorithms.}
\centering
	\begin{tabular}{ccccc}
		\toprule
		\multirow{2}{*}{N$\times$N} &
		\multicolumn{2}{c}{Computational time (s)} & \multirow{2}{*}{Speed-up ratio} & \multirow{2}{*}{$\norm{P_{FS}-P}_F$} \\
		
		\cline{2-3} &FS-2 &IPOT \\
		\midrule
		20$\times$20  & $2.24\times10^{-2}$ & $4.88\times10^{-1}$  & $2.18\times10^{1}$ &$2.40\times10^{-16}$   \\ 
		40$\times$40  & $7.74\times10^{-2}$ & $1.34\times10^{1}$ & $1.73\times10^{2}$ &$1.52\times10^{-16}$   \\ 
		80$\times$80  & $3.38\times10^{-1}$ & $4.79\times10^{2}$ & $1.42\times10^{3}$ &$1.40\times10^{-16}$   \\
		160$\times$160  & $1.42\times10^{0}$ & $1.46\times10^{4}$ & $1.03\times10^{4}$ &$8.51\times10^{-17}$   \\
		320$\times$320  & $6.31\times10^{0}$ & $-$ & $-$ & $-$   \\
		\bottomrule
	\end{tabular}
	\label{table:2drandom}
\end{table}

\subsection{Image matching problem}

\begin{figure} 
    \centering
    \includegraphics[width=0.6\linewidth]{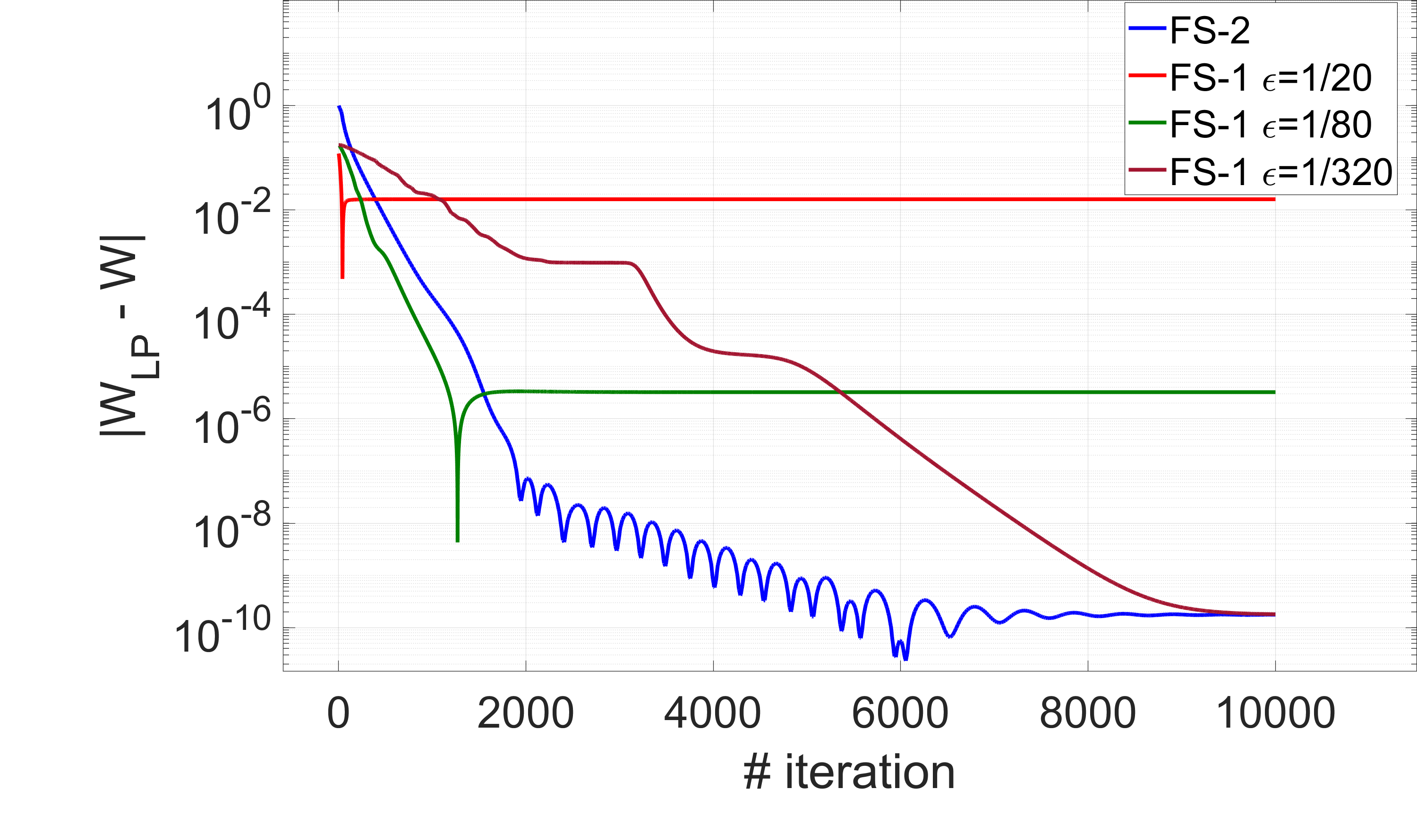}
    \caption{The image matching problem. The errors between the numerical results generated by FS-1 or FS-2 and the true Wasserstein-1 metric w.r.t. number of iterations.}
    \label{fig:graphiter}
\end{figure}

\begin{figure}
    \centering
    \includegraphics[width=0.5\linewidth]{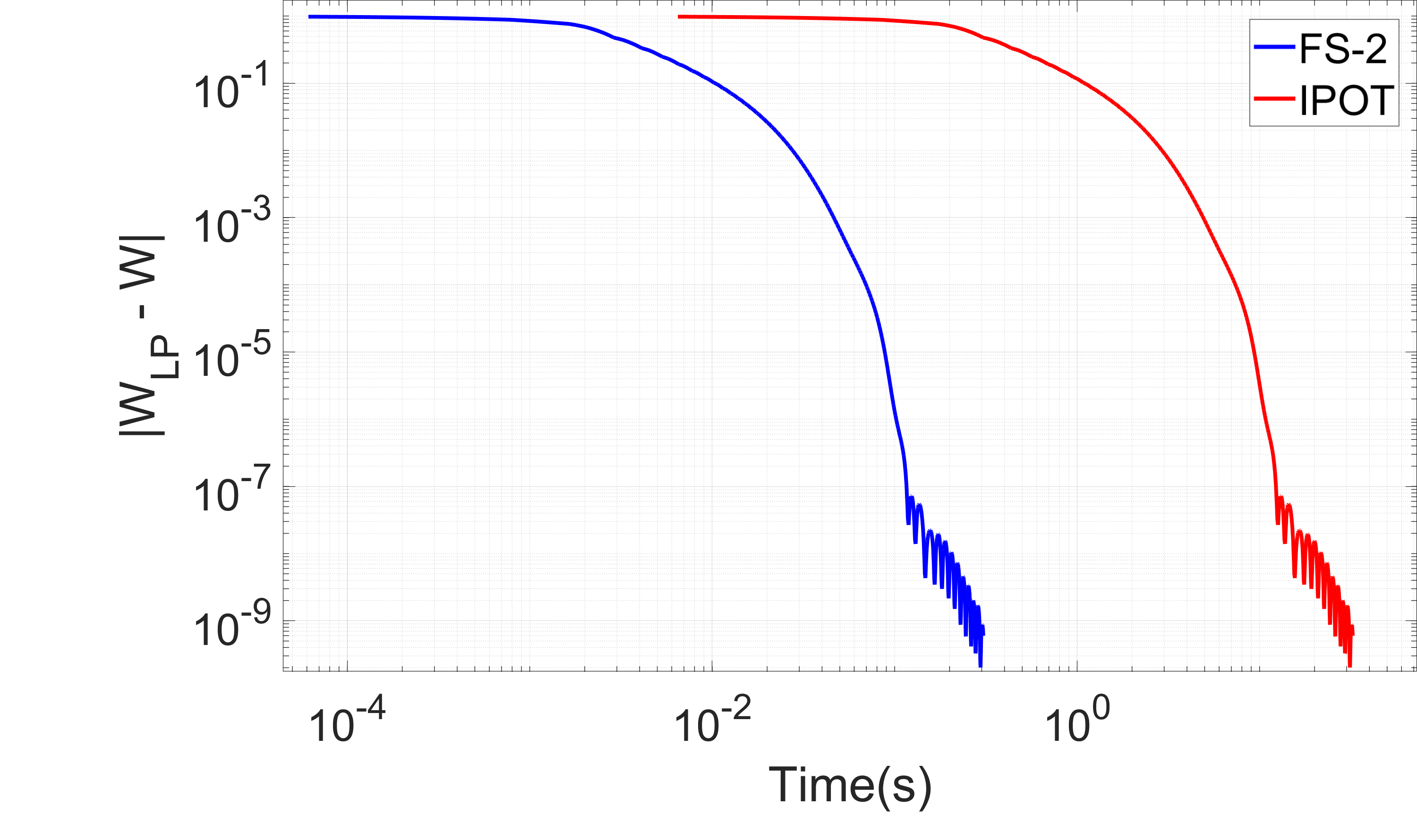}
    \caption{The image matching problem. The computational time required to reach the absolute error of the Wasserstein-1 metric.}
    \label{fig:2dfig}
\end{figure}

The final experiment tests the performance of our FS-2 algorithm for high-resolution image matching. This is a successful application of the Optimal Transport~\cite{burger2012regularized,haker2004optimal,clarysse2010optimal}. We select two images from the DIV2K dataset~\cite{Agustsson_2017_CVPR_Workshops}. Through a process similar to Subsection 5.4 in the manuscript~\cite{liao2022fast}, we compute the Wasserstein-1 metric between the two images.  The differences in the Wasserstein-1 metric between the true solution $W_{LP}$ and the numerical solutions generated by FS-1 and FS-2 are depicted in \cref{fig:graphiter}. We also present the averaged computational time of the IPOT method and the FS-2 algorithm in \cref{tab:2dfigtable}. Moreover, the computational time required to reach the absolute error of the Wasserstein-1 metric for $N\times N=32 \times 32$ is shown in \cref{fig:2dfig}. From these results, we can get the same conclusion as before, that is, the FS-2 algorithm seems to be the numerical algorithm with the fastest convergence and the lowest complexity for computing the Wasserstein-1 metric.

\begin{figure}
\centering
	\includegraphics[width=0.3\linewidth,height=0.25\linewidth]{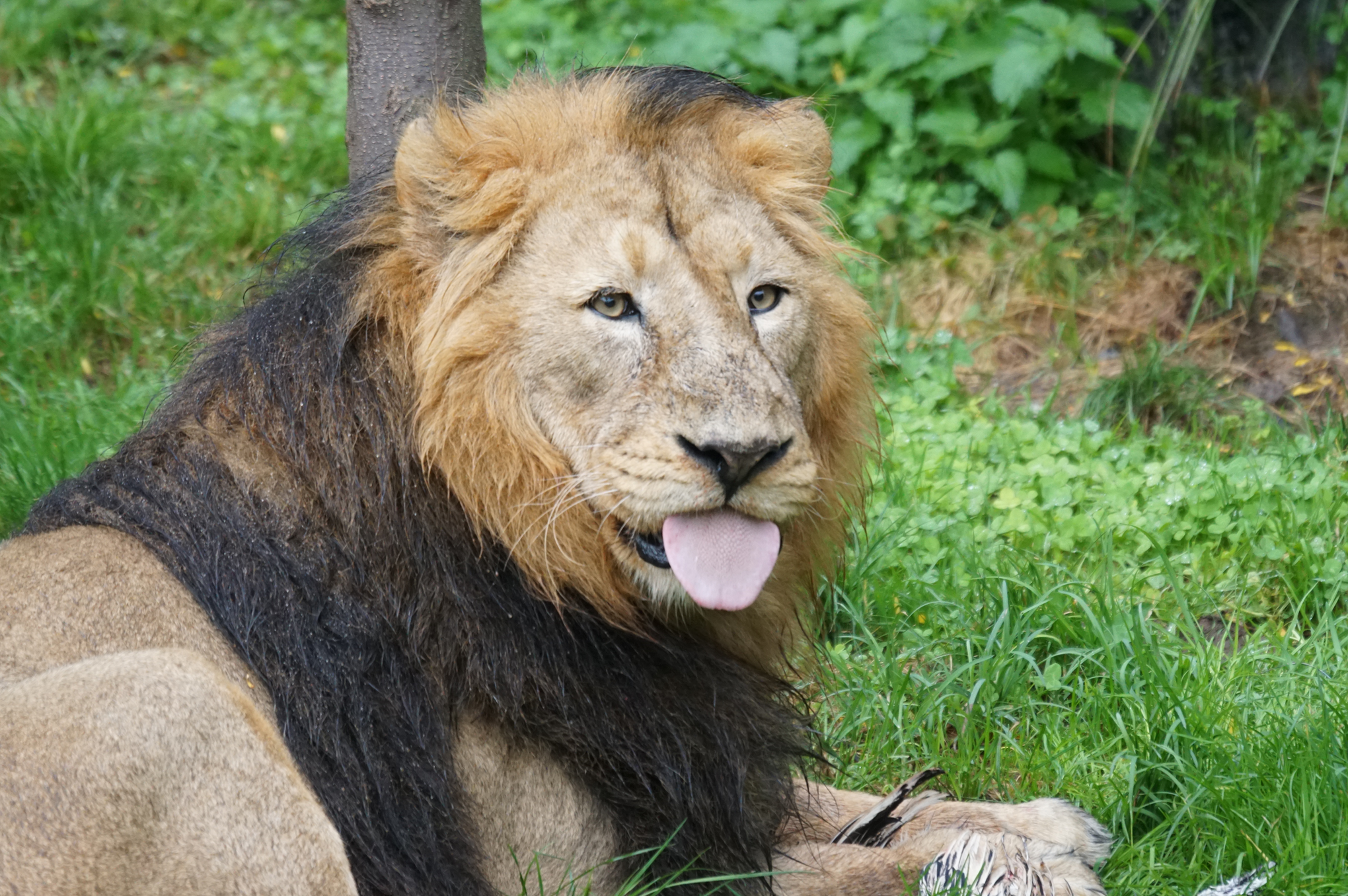}
	\includegraphics[width=0.3\linewidth,height=0.25\linewidth]{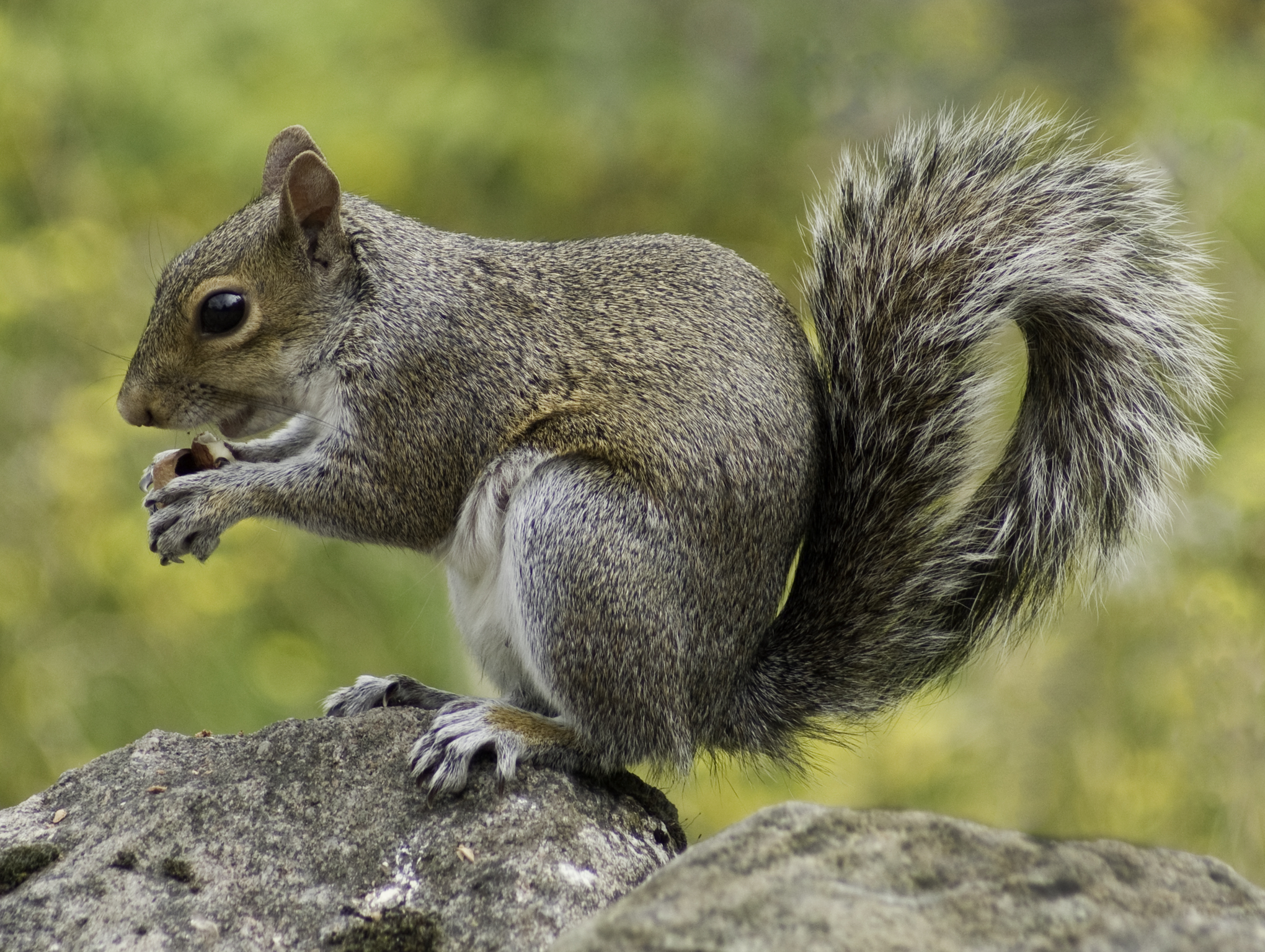}
\caption{The image matching problem. Illustration of images.}
\label{fig:ilsvrc}
\end{figure}

\begin{table}
	\caption{The image matching problem. The comparison between the IPOT method and the FS-2 algorithm with the different total number of grid nodes $N\times N$. Columns 2-4 are the averaged computational time of the two algorithms and the speed-up ratio of the FS-2 algorithm. Column 5 is the Frobenius norm of the difference between the transport plan computed by the two algorithms.}
\centering
	\begin{tabular}{ccccc}
		\toprule
		\multirow{2}{*}{$N\times N$} &
		\multicolumn{2}{c}{Computational time (s)} & \multirow{2}{*}{Speed-up ratio} & \multirow{2}{*}{$\norm{P_{FS}-P}_F$} \\
		
		\cline{2-3} & FS-2 &	IPOT \\
		
		 \midrule
		100$\times$100  & $5.16\times10^{-1}$ & $1.44\times10^{3}$ & $2.79\times10^{3}$  & $6.44\times10^{-17}$   \\ 
		200$\times$200  & $2.31\times10^{0}$ & $3.69\times10^{4}$  & $1.60\times10^{4}$ &$4.65\times10^{-17}$   \\ 
		400$\times$400  & $9.69\times10^{0}$ & $-$ & $-$ &$-$   \\ 
		800$\times$800  & $4.18\times10^{1}$ & $-$ & $-$ &$-$   \\
		\bottomrule
	\end{tabular}
	\label{tab:2dfigtable}
\end{table}

\section{Conclusion}
\label{sec:6}
As the follow-up of the FS-1 paper, we generalize the result of matrix-vector multiplication at $O(N)$ costs for the special matrix to the more general L/U-CoLT matrix. We illustrate that only two vectors are required to represent any L/U-CoLT matrix. Moreover, we also prove the closure of families of L/U-CoLT matrices to matrix Hadamard product and matrix scaling. Therefore, the above matrix operations are essentially updating the representation vectors, which reduce both time and space complexity to $O(N)$. These results can be directly applied to the Inexact Proximal point method for Optimal Transport problem and reduce the overall computational complexity to O(N). From this, we develop the Fast Sinkhorn II algorithm. It does not seem to be an overstatement that for the computation of the Wasserstein-1 metric, we have probably obtained the most competitive method, both in terms of convergence speed and computational complexity. 




\bibliographystyle{siamplain}
\bibliography{references.bib}
\end{sloppypar}
\end{document}

%% file: ex_shared.tex
\usepackage{lipsum}
\usepackage{amsfonts}
\usepackage{graphicx}
\usepackage{epstopdf}
\usepackage[noend]{algpseudocode}
\usepackage{algorithm,algorithmicx}
\usepackage{booktabs,multirow}
\ifpdf
  \DeclareGraphicsExtensions{.eps,.pdf,.png,.jpg}
\else
  \DeclareGraphicsExtensions{.eps}
\fi

\newsiamremark{remark}{Remark}
\newsiamremark{hypothesis}{Hypothesis}
\crefname{hypothesis}{Hypothesis}{Hypotheses}
\newsiamthm{claim}{Claim}

\headers{THE COLLINEAR TRIANGULAR MATRIX AND FS-2 ALGORITHM}{Q. LIAO, Z. WANG, J. CHEN, B. BAI, S. JIN AND H. WU}

\title{Fast Sinkhorn II: Collinear Triangular Matrix and Linear Time Accurate Computation of  Optimal Transport\thanks{Submitted to the editors DATE.
\funding{This work was supported by National Natural Science Foundation of People’s Republic of China Grant Nos.11871297 and 12031013, Tsinghua University Initiative Scientific Research Program, and Shanghai Municipal Science and Technology Major Project 2021SHZDZX0102.}}}

\author{Qichen Liao\thanks{Department of Mathematical Sciences, Tsinghua University, Beijing 100084, People’s Republic of China; Theory Lab, Central Research Institute, 2012 Labs, Huawei Technologies Co. Ltd., Hong Kong SAR, People’s Republic of China (\email{lqc20@mails.tsinghua.edu.cn}).}
\and Zihao Wang\thanks{Department of Computer Science and Engineering, Hong Kong University of Science and Technology, Clear Water Bay, Hong Kong SAR, People’s Republic of China (\email{zwanggc@cse.ust.hk}).}
\and Jing Chen\thanks{School of Physical and Mathematical Sciences, Nanyang Technological University, Singapore 639798 (\email{jing.chen@ntu.edu.sg}).}
\and Bo Bai\thanks{Theory Lab, Central Research Institute, 2012 Labs, Huawei Technologies Co. Ltd., Hong Kong SAR, People’s Republic of China (\email{baibo8@huawei.com}).}
\and Shi Jin\thanks{School of Mathematical Sciences, Institute of Natural Sciences, and MOE-LSC Shanghai Jiao Tong University, Shanghai 200240, People’s Republic of China (\email{shijin-m@sjtu.edu.cn}).}
\and Hao Wu\thanks{Corresponding author. Department of Mathematical Sciences, Tsinghua University, Beijing 100084, People’s Republic of China (\email{hwu@tsinghua.edu.cn})}.}

\usepackage{amsopn}
